\documentclass[english,russian]{amsart}
\usepackage{amssymb,amsmath}

\newtheorem{thm}{Theorem}[section]

\newtheorem{lm}[thm]{Lemma}

\newcommand{\der}{\partial}

\newcommand{\sudda}[1]{}

\begin{document}

\title{Asssociative algebras under multi-commutators}

\author{A.S. Dzhumadil'daev}

\address
{Kazakh-British University, Tole bi, 59, Almaty,  050000,
Kazakhstan} \email{dzhuma@math.kz}



\maketitle

\bigskip

\begin{abstract}  For an associative algebra $A$
a skew-symmetric (symmetric) sum of $n!$ products of $n$ elements of $A$ in all possible order is called Lie (Jordan) $n$-commutator. We consider $A$ as $n$-ary algebra under $n$-commutator.  
We construct $n$-ary skew-symmetric and symmetric generalizations of Jordan identity. 
We prove that any associative algebra under Jordan $n$-commutator satisfies a symmetric generalization of Jordan identity. We prove that in case of odd $n$  any associative algebra  
under Lie $n$-commutator satisfies a skew-symmetric generalization of Jordan identity. 
In case of even $n$ Lie $n$-commutator satisfies the homotopical $n$-Lie identity. 
 \end{abstract}

Well known that an associative algebra $A$ under Lie commutator is Lie. In other words, 
a vector space $A$ under commutator $[a,b]=ab-ba$ has skew-symmetric multiplication  $[\;,\;]:\wedge^2A\rightarrow A,$  that satisfies the identity, called Jacobi identity
$$[a_1,[a_2,a_3]]-[a_2,[a_1,a_3]]+[a_3,[a_1,a_2]]=0.$$ 

Well known also, that 
an associative algebra $A$ under Jordan  commutator $\{a,b\}=ab+ba$ is Jordan.  In other words, 
Jordan commutator is symmetric multiplication  $\{\;,\;\}:S^2A\rightarrow A,$  that satisfies the identity of degree 4, called Jordan identity
$$ \{a_1,\{a_0,\{a_2,a_3\}\}\}+\{a_2,\{a_0,\{a_1,a_3\}\}\}+\{a_3,\{a_0,\{a_1,a_2\}\}\}$$
$$-\{\{a_0,a_1\},\{a_2,a_3\}\}-\{\{a_0,a_2\},\{a_1,a_3\}\}-\{\{a_0,a_3\},\{a_1,a_2\}\}=0.$$

In our paper we consider multi-versions of these connections. We answer to a question of 
A.G. Kurosh who asked about identites of multi-associative algebras under multi-commutator \cite{Kurosh}. We show that an associative algebra under skew-symmetric $n$-commutator satisfies a  
homotopy identity (generalisation of Jacobi identity) if $n$ is even and one skew-symmetric generalization of Jordan identity  if $n$ is odd. We establish that an associative algebra under symmetric $n$-commutator satisfies symmetric generalization of Jordan identity. 

To formulate our results we need to introduce some definitions. 
Let $A$ be a vector space over a field $K.$ For a multilinear map $\alpha: A\times \cdots \times A\rightarrow A$ we say that  $A=(A,\alpha)$ is {\it $n$-algebra} with  {\it $n$-multiplication} $\alpha.$ 
A $n$-algebra $A$ is said {\it skew-commutative} if $\alpha$ is skew-symmetric, 
$$\alpha(a_{\sigma(1)},\ldots,a_{\sigma(n)})= sign\,\sigma\, \alpha(a_1,\ldots,a_n),$$
for any permutation $\sigma\in Sym_n.$ Similarly, $(A,\alpha)$ is {\it commutative} $n$-algebra, if 
$$\alpha(a_{\sigma(1)},\ldots,a_{\sigma(n)})= \alpha(a_1,\ldots,a_n),$$
for any $\sigma\in Sym_n.$

An absolute free $n$-algebra (free $n$-magma) can be defined as algebra of 
($n$-non-commutative, $n$-non-associative) $n$-polynomials  $K\langle t_1,t_2,\ldots \rangle.$ 
Denote by $\omega$ a $n$-multiplication in free $n$-magma. 
To construct  $n$-polynomials we have to introduce $n$-monoms. 

By definition, any variable $t_i$ is a $n$-monom of {\it $\omega$-degree $0.$ }
If $f_i$ is a $n$-monom of {\it $\omega$-degree $k_i,$} and $i=1,\ldots,N,$ then 
$\omega(f_1,\ldots,f_N)$ is a $n$-monom of {\it $\omega$-degree $k_1+\cdots+k_N+1.$} 
A linear combination of $n$-monoms is called ($n$-non-commutative, $n$-non-associative) $n$-polynomial. A space of $n$-polynomials  $K\langle t_1,t_2,\ldots \rangle$  is defined as a linear space with base generated by $n$-monoms. A multiplication $\omega$ on  
$K\langle t_1,t_2,\ldots \rangle$ is defined in a natural way. If $g_1,\ldots,g_N\in K\langle t_1,t_2,\ldots \rangle,$ then by multilinearity $\omega(g_1,\ldots,g_N)$ is a linear combination of $n$-monoms. We can imagine $n$-monoms as a rooted tree, where each vertex has $n$-in edges and $1$-out edge. Leaves are labeled by  elements of algebra and  to innner vertices  correspond $n$-ary products of elements that come by in-edges. 

Let  $f=f(t_1,\ldots,t_k) $ be any $n$-polynomial of $K\langle t_1,t_2,\ldots \rangle.$
Let $(A,\alpha)$ be any $n$-algebra with $n$-multiplication $\alpha.$ For any $k$ elements 
$a_1,\ldots,a_k\in A$ one can make substitutions $t_i:=a_i$ and $\omega:=\alpha$  in polynomial $f$ and consider another element $f(a_1,\ldots,a_k)$ of $A$ where multiplications  are made in terms of multiplication $\alpha$  instead of $\omega.$ We say that $f=0$ is a {\it  $n$-identity}  on 
$A$ if $f(a_1,\ldots,a_k)=0$ for any $a_1,\ldots,a_k\in A.$ 

In case of $n=2$ we obtain usual algebras.  In $2$-algebras multiplications are usually denoted as $a\circ b$ $a\times b$ $a+b,$ etc,  instead of $\alpha(a,b).$ The notions of $2$-polynomials and $2$-polynomial identities  are coincide with usual notions of polynomials and  polynomial identities

Let $$s_n=s_n(t_1,\ldots,t_n)=\sum_{\sigma\in Sym_n} sign\,\sigma\, t_{\sigma(1)}\cdots t_{\sigma(n)}$$
be a standard (associative, non-commutative) skew-symmetric polynomial. Then any associative algebra $A$ with  $2$-multiplication $ab$ can be endowed by a structure of $n$-algebra given by $n$-multiplication $s_n(a_1,\ldots,a_n).$ Call 
$$[a_1,\ldots,a_n]=s_n(a_1,\ldots,a_n)$$ as {\it Lie  $n$-commutator.} Note that Lie $2$-commutator coincides with usual Lie commutator,
$$s_2(a_1,a_2)=a_1a_2-a_2 a_1.$$

Let $$s_n^+=s_n^+(t_1,\ldots,t_n)=\sum_{\sigma\in Sym_n}  t_{\sigma(1)}\cdots t_{\sigma(n)}$$
be a standard (associative, non-commutative) symmetric polynomial. We can  endow  any associative algebra $A$ with  $2$-multiplication $ab$ by a structure of $n$-algebra given by  
{\it Jordan $n$-commutator}
$$\{a_1,\ldots,a_n\}=s_n^+(a_1,\ldots,a_n).$$  Note that Jordan $2$-commutator coincides with usual Jordan commutator,
$$s_2^+(a_1,a_2)=a_1a_2+a_2a_1.$$

In our paper we study $n$-polynomial identities of the algebra 
$(A,s_n),$ and $(A,s_n^+)$ if $2$-algebra $A$ is associative. In fact we construct generalizations of Lie and Jordan identities that hold for total associative algebras under Lie and Jordan $n$-commutators. 

\section{\bf Formulations of main results}

Let $Sym_{n}$ be set of all permutations on $[n]=\{1,2,\ldots,n\}.$ Let 
$$S_{k,l}=\{\sigma\in Sym_{k+l} | \sigma(1)<\ldots<\sigma(k), \sigma(k+1)<\cdots<\sigma(k+l)\}$$

$$S_{n-1,n-1,n}=\{\sigma\in Sym_{3n-2} | \sigma(1)<\ldots<\sigma(n-1),$$
$$ \sigma(n)<\cdots<\sigma(2n-2),\quad  \sigma(2n-1)<\cdots<\sigma(3n-2)\},$$

$$S_{n-2,n,n}=\{\sigma\in Sym_{3n-2} | 
 \sigma(1)<\ldots<\sigma(n-2),$$ $$ \sigma(n-1)<\cdots<\sigma(2n-2), \quad \sigma(2n-1)<\cdots<\sigma(3n-2), \quad \sigma(n-1)<\sigma(2n-1)\}.$$
are subsets of shuffle-permutations

If $n$-multiplication $\omega(t_1,\ldots,t_n)$ is skew-symmetric, then there exists  only one  
$n$-monomial of $\omega$-degree $2$
$$H(t_1,\ldots,t_{2n-1})=\omega(t_1,\ldots,t_{n-1},\omega(t_n,\ldots,t_{2n-1})).$$
Let 
$$h(t_1,\ldots,t_{2n-1})=\frac{1}{(n-1)!n!}\, \omega(t_{[1},\ldots,t_{n-1},\omega(t_{n},\ldots,t_{(2n-1)]}))$$
be its skew-symmetrisation by all parameters.  Note that,
$$h(t_1,\ldots,t_{2n-1})= \sum_{\sigma\in Sym_{n-1,n}} sign\,\sigma\,\omega(t_{\sigma(1)},\ldots,t_{\sigma(n-1)},(t_{\sigma(n)},\ldots,t_{\sigma(2n-1)})).$$

For skew-commutative $n$-algebras there are two $n$-monomials of $\omega$-degree $3$
$$F_1(t_1,\ldots,t_{3n-2})=(t_1,t_{2},\ldots,t_{n-1},(t_n,\ldots,t_{2n-2},(t_{2n-1},\ldots,t_{(3n-2)})))$$ and 
$$F_2(t_1,\ldots,t_{3n-2})=(t_1,t_{2},\ldots,t_{n-2},(t_{n-1},\ldots,t_{2n-2}),(t_{2n-1},\ldots,t_{(3n-2)})))
$$
Let us introduce their skew-symmetric sums by all parameters except $t_1,$  
$$F^{[2]}_1(t_1,\ldots,t_{3n-2})=\frac{1}{(n-1)!(n-2)!n!}\omega(t_{[2},t_3,\ldots,t_{n},\omega(t_{[1]},t_{n+1},\ldots,t_{2n-2},\omega(t_{2n-1},\ldots,t_{(3n-2)]}))),$$
$$
F^{[2]}_2(t_1,\ldots,t_{3n-2})=\frac{1}{(n-2)!(n-1)!n!}\omega(t_{[2},t_3,\ldots,t_{n-1},\omega(t_{[1]},t_{n},\ldots,t_{2n-2}),\omega(t_{2n-1},\ldots,t_{(3n-2)]})).
$$
Upper index $s$ in $F^{[s]}_l$ corresponds to the $\omega$-place where $t_1$ is  and lower index $l$ corresponds to $n$-bracketing types of $\omega$-degree $3.$
We have, 
$$F^{[2]}_1(t_1,\ldots,t_{3n-2})=$$
$$\sum_{\sigma\in Sym_{n-1,n-1,n}, \sigma(n)=1}sign\,\sigma\,
(t_{\sigma(1)},\ldots,t_{\sigma(n-1)},(t_{1},t_{\sigma(n+1)},\ldots,t_{\sigma(2n-2)},(t_{\sigma(2n-1)},\ldots,t_{\sigma(3n-2)}))),$$

$$F^{[2]}_2(t_1,\ldots,t_{3n-2})=$$
$$\sum_{\sigma\in Sym_{n-2,n,n}, \sigma(n-1)=1}sign\,\sigma\,
(t_{\sigma(1)},\ldots,t_{\sigma(n-2)},(t_{1},t_{\sigma(n)},\ldots,t_{\sigma(2n-2)}),(t_{\sigma(2n-1)},\ldots,t_{\sigma(3n-2)}))).$$
Let 
$$f_{\lambda}^{[2]}=F^{[2]}_1+\lambda F^{[2]}_2.$$

These notions have symmetric analoges. We save the same notations as in skew-symmetric case.  Just change brackets of the form $[ \,\;,\;]$ to $\{\, \;,\;\}.$ 

For commutative $n$-algebras there are two $n$-monomials of $\omega$-degree $3$
$$F_1^+(t_1,\ldots,t_{3n-2})=(t_1,t_{2},\ldots,t_{n-1},(t_n,\ldots,t_{2n-2},(t_{2n-1},\ldots,t_{(3n-2)})))$$ and 
$$F_2^+(t_1,\ldots,t_{3n-2})=(t_1,t_{2},\ldots,t_{n-2},(t_{n-1},\ldots,t_{2n-2}),(t_{2n-1},\ldots,t_{(3n-2)})))
$$
Their symmetric sums by all parameters except $t_1$ are defined by   
$$F^{\{2\}}_1(t_1,\ldots,t_{3n-2})=\frac{1}{(n-1)!(n-2)!n!}\omega(t_{\{2},t_3,\ldots,t_{n},\omega(t_{\{1\}},t_{n+1},\ldots,t_{2n-2},\omega(t_{2n-1},\ldots,t_{(3n-2)\}}))),$$
$$
F^{\{2\}}_2(t_1,\ldots,t_{3n-2})=\frac{1}{(n-2)!(n-1)!n!}\omega(t_{\{2},t_3,\ldots,t_{n-1},\omega(t_{\{1\}},t_{n},\ldots,t_{2n-2}),\omega(t_{2n-1},\ldots,t_{(3n-2)\}})).
$$
Upper index $s$ in $F^{\{s\}}_l$ corresponds to the place of $\omega,$ where $t_1$ is  and lower index $l$ corresponds to $n$-bracketing types of $\omega$-degree $3.$
We have, 
$$F^{\{2\}}_1(t_1,\ldots,t_{3n-2})=$$
$$\sum_{\sigma\in Sym_{n-1,n-1,n}, \sigma(n)=1}
(t_{\sigma(1)},\ldots,t_{\sigma(n-1)},(t_{1},t_{\sigma(n+1)},\ldots,t_{\sigma(2n-2)},(t_{\sigma(2n-1)},\ldots,t_{\sigma(3n-2)}))),$$

$$F^{\{2\}}_2(t_1,\ldots,t_{3n-2})=$$
$$\sum_{\sigma\in Sym_{n-2,n,n}, \sigma(n-1)=1}
(t_{\sigma(1)},\ldots,t_{\sigma(n-2)},(t_{1},t_{\sigma(n)},\ldots,t_{\sigma(2n-2)}),(t_{\sigma(2n-1)},\ldots,t_{\sigma(3n-2)}))).$$
Let 
$$f_{\lambda}^{\{2\}}=F^{\{2\}}_1+\lambda F^{\{2\}}_2.$$

Let $A$ be  $n$-algebra with $n$-multiplication $(a_1,\ldots,a_n).$ Denote by $[A]$ an algebra with vector space $A$ and $n$-multiplication  
$$[a_1,\ldots,a_n]=(a_{[1},\ldots,a_{n]})=\sum_{\sigma\in Sym_n}sign\,\sigma\,(a_{\sigma(1)},\ldots,a_{\sigma(n)})$$
(Lie $n$-commutator).  Similarly, denote by $\{A\}$ an algebra with vector space $A$ and $n$-multiplication  
$$\{a_1,\ldots,a_n\}=(a_{\{1},\ldots,a_{n\}})=\sum_{\sigma\in Sym_n}(a_{\sigma(1)},\ldots,a_{\sigma(n)})$$
(Jordan $n$-commutator). 

Recall that $A$ is called {\it total associative} \cite{Gnedbaye} if 
$$(a_1,\ldots,a_i,(a_{i+1},a_{i+2},\ldots,a_{i+n}),a_{i+n+1},a_{i+n+2},\ldots,a_{2n-1})=$$ $$
(a_1,\ldots,a_i,a_{i+1},(a_{i+2},\ldots,a_{i+n},a_{i+n+1}),a_{i+n+2},\ldots,a_{2n-1}),$$
for any $1\le i\le n-2.$
Any associative algebra $A$ under $n$-multiplication $(a_1,\ldots,a_n)\mapsto a1\cdot \cdots \cdot a_n$ became total associative.

The following skew-symmetric $\omega$-degree 2 polynomial is called {\it homotopical $n$-Lie}
$$homot(t_1,\ldots,t_{2n-1})=\sum_{\sigma\in S_{n-1,n}} sign\,\sigma\,\omega(t_{\sigma(1)},
\ldots, t_{\sigma(n-1)},\omega(t_{\sigma(n)},\ldots,t_{\omega(2n-1)})).$$
An $n$-ary algebra $(A,\omega)$ is called {\it homotopical $n$-Lie}, if it satisfies the identity $homot=0$ \cite{Hanlon}.

\begin{thm} \label{1} Let $A$ be total associative $n$-algebra. If $n$ is even or if 
$char \,K=p>0$ and $n$ is divisable by $p,$ then 
 its {$n-\!\!$commutators} algebra $[A]$ is homotopical $n$-Lie.
\end{thm}

\begin{thm} \label{2} Let $A$ be total associative algebra. Then its  Lie $n$-commutators algebra $[A]$ satisfies the identity 
$f_{-1}^{[2]}=0.$
\end{thm}

\begin{thm} \label{3} Let $A$ be total associative algebra. Then its Jordan $n$-commutators algebra $\{A\}$ satisfies the identity 
$f_{-1}^{\{2\}}=0.$
\end{thm}

{\bf Remarks.} The fact that $3$-commutatiors algebra $[A]$ has no identity of $\omega$-degree $2$ was noticed by A.G. Kurosh in \cite{Kurosh}. The identity $f_{-1}^{[2]}=$ holds for any $n$-commutators algebra, but this identity in general is not minimal. If $n$ is even or if the characteristic of main field is $p>0$ and $n\equiv  0(mod\,p),$  then one can find for $[A]$  the identity of $\omega$-degree 2, for example,  $homot=0.$ 
We think   that $homot=0$  for even $n$ and $f^{[2]}_{-1}=0$ for odd $n$ are minimal identities that hold  for any  Lie $n$-commutator algebras $[A],$  if $char K=0.$ We think also that $f_{-1}^{\{2\}=0}$ is minimal identity that hold for any Jordan $n$-commutators algebra $\{A\}.$

The case $n=3$ was considered by M. R. Bremner \cite{Bitc}, \cite{BPtalma}.
He proved that $f^{[2]}_{-1}=0$ and $f^{\{2\}}_{-1}=0$ are identities   for Lie and Jordan 
3-commutators of  total associative algebras and 
he established  the minimality of these identities. 

In case of $n=2$  the polynomial $f_{-1}^{\{2\}}$ coincides with usual  Jordan polynomial.

\section{\bf Proof of Theorem \ref{1}}

Let  $A$ be a free total associative $n$ algebra with $n$-multiplication $\omega$ and $[\omega]$ be its $n$-commutator,
$$[\omega](t_1,\ldots,t_n)=\sum_{\sigma\in Sym_n} sign\,\sigma\, \omega(t_{\sigma(1)},\cdots, t_{\sigma(n)})$$
We have to prove that $X=0,$ 
where
$$X=X(t_1,\ldots,t_{2n-1})=\sum_{\sigma\in S_{n-1,n}}
sign\,\sigma\,[\omega](t_{\sigma(1)},\ldots,t_{\sigma(n-1)},[\omega](t_{\sigma(n)},\ldots,t_{\sigma(2n-1)})).$$
Expand  $n$-commutators $[\omega]$ in terms of associative $n$-multiplication $\omega$. We see 
 that $X$ is a sum of elements of a form 
$$\pm \omega(t_{i_1},\ldots,t_{i_s},\omega(t_{i_{s+1}},\ldots,t_{i_{s+n}}),t_{i_{s+n+1}},\ldots,t_{2n-1}).$$
Since $A$ is total associative, this sum is reduced to a  sum of elements of a form 
$$\pm \omega(t_{j_1},\ldots,t_{j_{n-1}},\omega(t_{j_n},\ldots,t_{j_{2n-1}})).$$
Let $\mu \in K$ be the coefficient of $X$ at $\omega(t_1,\ldots,t_{n-1},\omega(t_n,\ldots,t_{2n-1}))$ 
Since $X(t_1,\ldots,t_{2n-1})$ is skew symmetric by all arguments $t_1,\ldots,t_{2n-1}$  
to prove $X=0$ it is enough to establish that $\mu=0.$

Note that the element $Q:=\omega(t_1,\ldots, t_{n-1},\omega(t_{n},\ldots,t_{2n-1}))$ 
may enter with non-zero coefficient only in summands of $X$ of a form 
$$R_{n-1}:=[\omega](t_1,\ldots,t_{n-1},[\omega](t_n,\ldots,t_{2n-1})),$$ $$
R_{n-2}:=[\omega](t_1,\ldots,t_{n-2},t_{2n-1},[\omega](t_{n-1},\ldots,t_{2n-2})),\quad \ldots,$$ $$
R_0:=[\omega](t_{n+1},\ldots,t_{2n-1},[\omega](t_1,\ldots,t_n)).$$
The element  
$$R_i=[\omega](t_1,\ldots,t_i,t_{n+i+1},\ldots t_{2n-1},[\omega](t_{i+1},\cdots,  t_{i+n})), \quad 0\le i\le n-1,$$
enter to $X$ with coefficient that is equal to signature of the permutation 
$$\gamma_i=\left(\begin{array}{ccccccccc}
1&\cdots&i&i+1&\cdots&{n-1}&n&\cdots&2n-1\\
1&\cdots&i&n+i+1&\cdots&2n-1&i+1&\cdots&i+n
\end{array}\right)\in S_{n-1,n}
$$
We have 
$$sign\,\gamma_i=(-1)^{(n-i-1)n}.$$
To obtain a component $Q$ from $R_i$ we have to permute the part 
$\omega(t_{i+1},\ldots,t_{i+n})$ of $R_i$  $(n-i-1)$ times,
$$[\omega](t_1,\ldots,t_i,t_{n+i+1},\ldots t_{2n-1},[\omega](t_{i+1},\cdots,  t_{i+n})) 
\rightsquigarrow \cdots $$
$$
\rightsquigarrow(-1)^{n-1-i}\omega(t_{1},\ldots,t_{i},\omega(t_{i+1},\ldots,t_{i+n}),t_{i+n+1},\ldots,t_{2n-1}).$$
Therefore
$$\mu=\sum_{i=0}^{n-1}sign\,\gamma_i(-1)^{n-i-1}=\sum_{i=0}^{n-1}(-1)^{(n-1-i)(n+1)}.$$
Note that 
$$\mu=\left\{\begin{array}{cc} 0&\mbox{ if $n$ is even}\\
n& \mbox{ if $n$ is odd }
\end{array}\right.$$
Hence, $X=0,$ if $n$ even or   $char K=p>0,$ $n$ is odd and  $n$ is divisable by $p.$

\section{\bf Proof of Theorem \ref{2}}

Note that 
$$F^{[2]}_1(t_1,\ldots,t_{3n-2})=
\sum_{\sigma\in Sym_{n-1,n-1,n}, \sigma(n)=1}$$
$$sign\,\sigma\,
[\omega](t_{\sigma(1)},\ldots,t_{\sigma(n-1)},[\omega](t_{1},t_{\sigma(n+1)},\ldots,t_{\sigma(2n-2)},[\omega](t_{\sigma(2n-1)},\ldots,t_{\sigma(3n-2)}))),$$

$$F^{[2]}_2(t_1,\ldots,t_{3n-2})=\sum_{\sigma\in Sym_{n-2,n,n}, \sigma(n-1)=1}$$ $$sign\,\sigma\,
[\omega](t_{\sigma(1)},\ldots,t_{\sigma(n-2)},[\omega](t_{1},t_{\sigma(n)},\ldots,t_{\sigma(2n-2)}),[\omega](t_{\sigma(2n-1)},\ldots,t_{\sigma(3n-2)}))).$$

For any permutation  $i_1i_2\ldots i_{3n-2}\in Sym_{3n-2}$ set  $$e(i_1\ldots i_{3n-2})
:=\omega(t_{i_1},\ldots,t_{i_{n-1}},\omega(t_{i_n},\ldots,t_{i_{2n-2}},\omega(t_{i_{2n-1}},\ldots,
t_{i_{3n-2}})))$$
and 
$$[e](i_1\ldots i_{3n-2})
:=[\omega](t_{i_1},\ldots,t_{i_{n-1}},[\omega](t_{i_n},\ldots,t_{i_{2n-2}},[\omega](t_{i_{2n-1}},\ldots,
t_{i_{3n-2}}))).$$
For any $1\le i\le 3n-2$ let 
$$e_i=e(2,\ldots,i,1,i+1,\ldots,3n-2)$$
The index $i$ corresponds the place where is $1$. For example, $e_1=e(1,2,\ldots,3n-~2), e_2=e(2,1,3,\ldots,3n-2), e_{3n-2}=e(2,\ldots,3n-2,1).$

Since $A$ is total associative,  for any $s=1,2,$ the element 
$F^{[2]}_s(t_1,\ldots,t_{3n-2})$  can be presented as a sum of elements 
$e(i_1,\ldots,i_{3n-2}),$ where $i_1\ldots i_{3n-2}$ is a  permutation of the set $[3n-2]=\{1,2,\ldots,3n-2\}.$ 

Let $\mu_{s,i}$ be a coefficient of $F^{[2]}_s(t_1,\ldots,t_{3n-2})$ at $e_i=
e(2,\ldots,i,1,i+1,\ldots,3n-2),$
where $1\le i\le 3n-2.$ 
Since  $F^{[2]}_s(t_1,t_2,\ldots,t_{3n-2})$ is skew-symmetric by all variables except $t_1,$ the element $F^{[2]}_s(t_1,t_2,\ldots,t_{3n-2})$ is uniquily defined by coefficient $\mu_{s,i},$ where $1\le i\le 3n-2.$  
Then the  condition that  $f^{[2]}_{-1}=0$ is identity on $[A]$ is equivalent to the following relations 
\begin{equation}\label{12}
\mu_{1,i}=\mu_{2,i},\quad 1\le i\le 3n-2.
\end{equation}
We will establish the following common values for  $\mu_{1,i}$ and $\mu_{2,i}.$

Let for even  $n$ 

$$\mu_{i}=\left\{ \begin{array}{ll}
(-1)^{(i+1)}\lfloor \frac{i+1}{2} \rfloor & \mbox{ if $i\le n$ }\\
&\\
(-1)^{i+1}(\frac{n}{2}+2\lfloor \frac{n-i-1}{2}\rfloor)& \mbox{ if $ n+1\le  i\le 2n-2$ }\\
&\\
(-1)^i \lfloor \frac{3n-i}{2}\rfloor& \mbox{ if $2n-1\le i\le 3n-2$ }\\
\end{array}\right.
$$
and for odd  $n$

$$\mu_{i}=\left\{ \begin{array}{ll}
(-1)^{(i+1)}\,\frac{(2n-i-1)i}{2} & \mbox{ if $i\le n$ }\\
&\\
(-1)^{(i+1)}\, (\frac{n(n-1)}{2}+(i-n)(-2n+i+1)) & \mbox{ if $ n+1\le  i\le 2n-2$ }\\
&\\
(-1)^i \, \frac{(3n-i-1)(n-i)}{2}& \mbox{ if $2n-1\le i\le 3n-2$ }\\
\end{array}\right.
$$
Note that 
$$\mu_{i}=\mu_{3n-1-i}, \qquad 1\le i\le 3n-2.$$

Let  $[i,j]=\{ s\in {\bf Z} | i\le s\le j\}$ be segment with endpoints $i,j$ and 
$[i,j)=\{s\in {\bf Z} | i\le s<j\},$ $(i,j]=\{s\in {\bf Z} | i<s\le j\},$ $(i,j)=\{s\in {\bf Z} | i<s<j\}$  be semi-segments. Note that semi-segment $[i,j)$ has endpoints $i$ and $j-1$ and similarly, endpoints of $(i,j]$ is $i+1$ and $j.$  
Number of elements  of (semi)-segment is called length. For example, 
$|[i,j]|=j-i$ and $|[i,j)|=j-1-i,$ if $j>i.$
Say that $[i_1,j_1]\subseteq [i,j]$ is subsegment if $i\le i_1<j_1\le j.$

\begin{lm}\label{40} 
Let $\mu_{1,i}$ be the coefficient at $e_i$ of the element 
$F^{[2]}_1(t_1,\ldots,t_{3n-2}).$ Then 
$$\mu_{1,i}=\mu_i,$$
for any $1\le i\le 3n-2.$ 
\end{lm}

{\bf Proof.}  Consider in the  segment $P_1=[2,3n-2]=\{2,\ldots,3n-2\}$ chain with two subsegements
$$P_3\subset P_2\subset P_1,  \qquad |P_1|=3n-3, |P_2|=2n-2, |P_3|=n$$
Denote endpoints of $P_1,P_2,P_3$ as $A_1,B_1$ $A_2,B_2$ and $A_3,B_3.$ 
Then 
$$P_1=[2,3n-2], P_2=[p+1,p+2n-2], P_3=[q,q+n-1]$$
for some $1\le p<q \le 2n-1$ and the points $A_1,A_2,A_3,B_3,B_2,B_1$ on ${\bf R}$  has 
coordinates $2,p+1,q,q+n-1,p+2n-2,3n-2.$  Note that 
\begin{equation}\label{border}
1\le p\le n,\quad  p< q\le 2n-1,\quad  q\le p+n-1
\end{equation}

Then 
$$P_1=[A_1,A_2)\cup P_2\cup (B_2,B_1],$$
$$P_2=[A_2,A_3)\cup P_3\cup (B_3,B_2].$$
Let us introduce the following subsets of increasing integers 
$$X_1=[A_1,A_2)\cup(B_2,B_1],$$
$$X_2=\{1\}\cup [A_2,A_3)\cup (B_3,B_2],$$
$$X_3=P_3=[A_3,B_3].$$
In the following picture parts of $X_1, X_2,X_3$ are marked equally. 

$\begin{minipage}{2cm}
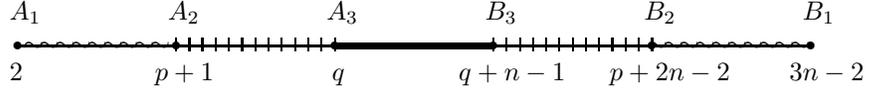
\begin{picture}(0,70)
\put(30,30){\line(1,0){300}}

\put(30,30){\circle*{3}}
\put(90,30){\circle*{3}}
\put(150,30){\circle*{3}}
\put(210,30){\circle*{3}}
\put(270,30){\circle*{3}}
\put(330,30){\circle*{3}}

\put(25,40){$A_1$}
\put(85,40){$A_2$}
\put(145,40){$A_3$}
\put(205,40){$B_3$}
\put(265,40){$B_2$}
\put(325,40){$B_1$}

\put(27,17){2}
\put(80,17){$p+1$}
\put(147,17){$q$}
\put(195,17){$q+n-1$}
\put(252,17){$p+2n-2$}
\put(320,17){$3n-2$}

\put(30,28){$\sim$}
\put(36,28){$\sim$}
\put(42,28){$\sim$}
\put(48,28){$\sim$}
\put(54,28){$\sim$}
\put(60,28){$\sim$}
\put(66,28){$\sim$}
\put(72,28){$\sim$}
\put(78,28){$\sim$}
\put(266,28){$\sim$}
\put(272,28){$\sim$}
\put(278,28){$\sim$}
\put(284,28){$\sim$}
\put(290,28){$\sim$}
\put(296,28){$\sim$}
\put(302,28){$\sim$}
\put(308,28){$\sim$}
\put(314,28){$\sim$}
\put(320,28){$\sim$}

\put(150,31){\line(1,0){60}}
\put(150,29){\line(1,0){60}}

\put(90,28){\line(0,1){5}}
\put(95,28){\line(0,1){5}}
\put(100,28){\line(0,1){5}}
\put(105,28){\line(0,1){5}}
\put(110,28){\line(0,1){5}}
\put(115,28){\line(0,1){5}}
\put(120,28){\line(0,1){5}}
\put(125,28){\line(0,1){5}}
\put(130,28){\line(0,1){5}}
\put(135,28){\line(0,1){5}}
\put(140,28){\line(0,1){5}}
\put(145,28){\line(0,1){5}}
\put(150,28){\line(0,1){5}}

\put(210,28){\line(0,1){5}}
\put(215,28){\line(0,1){5}}
\put(220,28){\line(0,1){5}}
\put(225,28){\line(0,1){5}}
\put(230,28){\line(0,1){5}}
\put(235,28){\line(0,1){5}}
\put(240,28){\line(0,1){5}}
\put(245,28){\line(0,1){5}}
\put(250,28){\line(0,1){5}}
\put(255,28){\line(0,1){5}}
\put(260,28){\line(0,1){5}}
\put(265,28){\line(0,1){5}}
\put(270,28){\line(0,1){5}}
\end{picture}
\end{minipage}
$

So, for any such chain $P_1\supset P_2\supset P_3$ one corresponds a sequence of elements 
$X_1X_2X_3$ where in each part $X_i$ elements are written in increasing order and $X_2$ begins by $1.$ In other words, any chain $P_1\supset P_2\supset P_3$ defines in a unique way an element $[e](X_1X_2X_3).$ More exactly, 
$$[e](X_1X_2X_3)=$$
$$[\omega](t_2,\ldots,t_{p},t_{p+2n-1},\ldots,t_{3n-2},[\omega](t_1,t_{p+1 },\ldots,t_{q-1},t_{q+n}\ldots,t_{p+2n-2},[\omega](t_q,\ldots,t_{q+n-1}))).$$
Signature of the permutation 
$$X_1X_2X_3 =
2\ldots p\,p+2n-1\,\ldots 3n-2\,1\,p+1\ldots q-1\,q+n\ldots p+2n-2\,q\ldots q+n-1$$
is equal to 
$$(-1)^{(p-1)+(n-p)(2n-1)+(n-q+p-1)n}.$$
So, 
\begin{equation}\label{13}
sgn\,X_1X_2X_3=(-1)^{(p+1-q)n+1}
\end{equation}

For $1\le i\le 3n-2$ and $1\le p\le n,$ $0<q-p\le n-1,$ denote by $\mu_{1,i}^{(p,q)},$
 the coefficient at  $e_i$
of the element  
$$[\omega](t_2,\ldots,t_{p},t_{p+2n-1},\ldots,t_{3n-2},[\omega](t_1,t_{p+1 },\ldots,t_{q-1},t_{q+n}\ldots,t_{p+2n-2},[\omega](t_q,\ldots,t_{q+n-1}))).$$
In case of $p=1,$ by $\mu_{1,i}^{(p,q)}$ we understand the coefficient at 
$e_i$ of the element 
$$[\omega](t_{2n},\ldots,t_{3n-2},[\omega](t_1,t_{2 },\ldots,t_{q-1},t_{q+n}\ldots,t_{p+2n-2},[\omega](t_q,\ldots,t_{q+n-1}))).$$
In case of $p=n,$ by $\mu_{1,i}^{(p,q)}$ we mean the coefficient at 
$e_i$ of the element 
$$[\omega](t_{2},\ldots,t_n,[\omega](t_1,t_{n+1},\ldots,t_{q-1},t_{q+n}\ldots,t_{p+2n-2},[\omega](t_q,\ldots,t_{q+n-1}))).$$
For any $1\le i \le 3n-2$ denote by $\mu_{1,i}^{(0,q)}$ the coefficient at $e_i$ of the element 
$$\omega(t_{2n},\ldots,t_{3n-2}, \omega(t_1,t_{n+2},\ldots,t_{2n-1},\omega(t_2,\ldots,t_{n+1}))).$$
Then $\mu_{1,i}^{(0,q)}=0,$ if $i\le n$ or $i\ge 2n-1.$

Below we use the following notation 
$Y\rightsquigarrow Z$ that means that $Z$ is obtained from $Y$ by using skew-symmetry property 
of $[\omega]$

We have, 
$$[\omega](t_2,\ldots,t_{p},t_{p+2n-1},\ldots,t_{3n-2},[\omega](t_1,t_{p+1},\ldots,t_{q-1},t_{q+n}\ldots,t_{p+2n-2},\omega(t_q,\ldots,t_{q+n-1})))\rightsquigarrow$$

$$(-1)^{(p-q+n-1)}
[\omega](t_2,\ldots,t_{p},t_{p+2n-1},\ldots,t_{3n-2},[\omega](t_1,
t_{p+1 },\ldots,t_{q-1},\omega(t_q,\ldots,t_{q+n-1}),t_{q+n}\ldots,t_{p+2n-2}))\rightsquigarrow$$

$$(-1)^{(q-1)}
[\omega](t_2,\ldots,t_{p},
[\omega](t_1,
t_{p+1 },\ldots,t_{q-1},\omega(t_q,\ldots,t_{q+n-1}),t_{q+n}\ldots,t_{p+2n-2}),
t_{p+2n-1},\ldots,t_{3n-2})\rightsquigarrow$$

$$(-1)^{(q-1)}
\omega(t_2,\ldots,t_{p},
[\omega](t_1,
t_{p+1 },\ldots,t_{q-1},\omega(t_q,\ldots,t_{q+n-1}),t_{q+n}\ldots,t_{p+2n-2}),
t_{p+2n-1},\ldots,t_{3n-2})$$
Now expand $[\omega]$  in 
$[\omega](t_1,
t_{p+1 },\ldots,t_{q-1},\omega(t_q,\ldots,t_{q+n-1}),t_{q+n}\ldots,t_{p+2n-2}).$  
Then $t_1$ might be in $i$-th place only in the following cases 

\begin{equation}\label{19}
\omega(
t_1,t_{p+1 },\ldots, t_{q-1},\omega(t_q,\ldots,t_{q+n-1}),t_{q+n}\ldots,t_{p+2n-2}),\quad i=p,
\end{equation}

\begin{equation}\label{20}
(-1)^{i-p}\omega(t_{p+1 },\ldots, t_i,t_1,\ldots, t_{q-1},\omega(t_q,\ldots,t_{q+n-1}),t_{q+n}\ldots,t_{p+2n-2}),\quad p+1\le i\le q-1,
\end{equation}

\sudda{\begin{equation}\label{200}
\omega(
t_{p+1 },\ldots, t_{q-1}, t_1,\omega(t_q,\ldots,t_{q+n-1}),t_{q+n}\ldots,t_{p+2n-2})
\end{equation}

\begin{equation}\label{201}
\omega(
t_{p+1 },\ldots, t_{q-1},\omega(t_q,\ldots,t_{q+n-1}),t_1,t_{q+n}\ldots,t_{p+2n-2}),\quad q+n-1\le i\le p+2n-2
\end{equation}
}

\begin{equation}\label{21} 
(-1)^{n-1-p+i}\omega(
t_{p+1 },\ldots,  t_{q-1},\omega(t_q,\ldots,t_{q+n-1}),t_{q+n}\ldots,t_i, t_1,\ldots, t_{p+2n-2}),\quad 
q+n-1\le i\le p+2n-2
\end{equation}

Note that $\mu_{1,i}^{(p,q)}=0,$ if $i\notin [A_2,A_3)\cup(B_3,B_2].$ Therefore,  by (\ref{19}), (\ref{20}) and (\ref{21}), 
\begin{equation}\label{14}
\mu_{1,i}^{(p,q)}=\left\{\begin{array}{ll}
0& \mbox{ if $i<p$ or $q\le i\le q+n-1$ or $p+2n-2\le i\le 3n-2$}\\
&\\
(-1)^{i+1 +p-q}& \mbox{ if $p\le i \le q-1$}\\  
&\\
(-1)^{n+p+i-q}& \mbox{ if $q+n-1\le i\le p+2n-2$}\\
\end{array}\right.
\end{equation}
Note also  $1\le p\le n, p< q\le p+n-1.$ Hence $q\le 2n-1.$

The element $e_i,$ where $1\le i\le 3n-2,$ may appear in expanding of 
$$[\omega](t_2,\ldots,t_{p},t_{p+2n-1},\ldots,t_{3n-2},[\omega](t_1,t_{p+1 },\ldots,t_{q-1},t_{q+n}\ldots,t_{p+2n-2},[\omega](t_q,\ldots,t_{q+n-1})))$$
with the coefficient 
\begin{equation}\label{15}\mu_{1,i}=\sum_{p,q} \mu_{1,i}^{(p,q)}.
\end{equation}

Let $i\le n.$ Then the case $q+n-1\le i,$ $1\le p<q$  is impossible. Therefore,
by (\ref{12}), (\ref{13}), (\ref{14}),
$$\mu_{1,i}=\sum_{p=1}^{i}\sum_{q=i+1}^{p+n-1}\mu_{1,i}^{(p,q)}
=\sum_{p=1}^i\sum_{q=i+1}^{p+n-1}(-1)^{(p+1-q)n+1}(-1)^{p-q+i+1}=
\sum_{p=1}^i\sum_{q=i+1}^{p+n-1}(-1)^{p-q+i+pn+qn+n}.$$
So, for  even  $n,$
$$\mu_{1,i}=
\sum_{p=1}^i\sum_{q=i+1}^{p+n-1}(-1)^{p-q+i}=
(-1)^{i}\sum_{p=1}^i \sum_{q=i+1}^{p+n-1} (-1)^{p-q}=(-1)^{i+1}\lfloor\frac{i+1}{2}\rfloor.
$$
For odd $n$
$$\mu_{1,i}=
\sum_{p=1}^i\sum_{q=i+1}^{p+n-1}(-1)^{i+1}=
(-1)^{i}\sum_{p=1}^i \sum_{q=i+1}^{p+n-1} 1=(-1)^{i+1}\frac{(2n-i-1)i}{2}.$$

Consider the case  $n+1\le i\le 2n-2.$ By  (\ref{12}), (\ref{13}), (\ref{14}),
$$\mu_{1,i}=\sum_{p=1}^n \sum_{q=i+1}^{p+n-1}(-1)^{(p+1-q)n+1}(-1)^{i+1+p-q}+\sum_{p=0}^n\sum_{q=p+1}^{i-n}(-1)^{(p+1-q)n+1}(-1)^{n+p+i-q}$$
So, if $n$ is even, then 
$$\mu_{1,i}=
\sum_{p=1}^n \sum_{q=i+1}^{p+n-1}(-1)^{i+p-q}-\sum_{p=0}^n\sum_{q=p+1}^{i-n}(-1)^{i-q+p}=
$$
$$(-1)^i(\sum_{p=1}^n \sum_{q=i+1}^{p+n-1}(-1)^{p-q}-\sum_{p=0}^n\sum_{q=p+1}^{i-n}(-1)^{q-p})=
$$
$$(-1)^{i+1}(\frac{n}{2}+2\lfloor \frac{n-i-1}{2}\rfloor).  $$
If $n$ is odd, then
$$\mu_{1,i}=\sum_{p=1}^n \sum_{q=i+1}^{p+n-1}(-1)^{i+1}+\sum_{p=0}^n\sum_{q=p+1}^{i-n}(-1)^{i-1}= 
$$
$$(-1)^{i + 1} (\frac{n (n - 1)}{2} + (i - n) (-2 n + i + 1)).$$

Consider the case $2n-1\le i\le 3n-2.$ Then all cases except  $q+n-1\le i\le p+2n-2$ are not possible. Therefore, $i-2n+2\le p<q\le i-n+1,$ and 
$$\mu_{1,i}=\sum_{p=i-2n+2}^n\sum_{q=p+1}^{i-n+1}\mu_{1,i}^{(p,q)}= 
\sum_{p=i-2n+2}^n \sum_{q=p+1}^{i-n+1}(-1)^{(p+1-q)n+1}(-1)^{n+p+i-q} $$
Hence, for even $n$
$$\mu_{1,i}=\sum_{p=i-2n+2}^n \sum_{q=p+1}^{i-n+1}(-1)^{i+p-q+1} =
(-1)^{i+1}\sum_{p=i-2n+2}^n \sum_{q=p+1}^{i-n+1}(-1)^{p-q} =
(-1)^i\lfloor \frac{3n-i}{2}\rfloor.$$
and for  odd $n$ 
$$\mu_{1,i}=\sum_{p=i-2n+2}^n\sum_{q=p+1}^{i-n+1}(-1)^{i+1}=
(-1)^{i+1}\sum_{p=i-2n+2}^n\sum_{q=p+1}^{i-n+1}1=
(-1)^i\frac{(3n-i-1)(n-i)}{2}.$$
Lemma \ref{40} is proved.

\begin{lm}\label{41} 
Let $\mu_{2,i}$ be the coefficient at $e_i$ of the element 
$F^{[2]}_2(t_1,\ldots,t_{3n-2}).$ Then 
$$\mu_{2,i}=\mu_i,$$
for any $1\le i\le 3n-2.$ 
\end{lm}

{\bf Proof.} 
Consider in the  segment $P_1=[2,3n-2]=\{2,\ldots,3n-2\}$ two non-intersecting subsegements 
of length $n-1$ and $n$
$$P_1\supset P_2,\quad  P_1\supset P_3,  \qquad |P_1|=3n-3, |P_2|=n-1  , |P_3|=n.$$
Denote endpoints of $P_1,P_2,P_3$ as $A_1,B_1$ $A_2,B_2$ and $A_3,B_3.$ 
Then 
$$P_1=[2,3n-2], \quad P_2=[p+1,p+n-1],\quad  P_3=[q,q+n-1]$$
for some $1\le p \le 2n-1$ and $2\le q\le 2n-1.$ 

Note that 
$$
q\ge p+ n \mbox{ if $p<q$}
$$
$$
p\ge q+ n-1 \mbox{ if $p>q$}
$$

Then 
$$P_1=[A_1,B_1],\quad P_2=[A_2,B_2],\quad  P_3=[A_3,B_3]$$

Let us introduce the following subsets of increasing integers 
$$X_1=[A_1,A_2)\cup(B_2,A_3)\cup (B_3,B_1]  \quad \mbox{\qquad (Case I)},$$
$$X_1=[A_1,A_3)\cup(B_3,A_2)\cup (B_2,B_1]  \quad \mbox{\qquad (Case II)},$$
$$X_2=\{1\}\cup [A_2,B_2],$$
$$X_3=P_3=[A_3,B_3].$$
In the following picture parts of $X_1, X_2,X_3$ are marked equally. 

$\begin{minipage}{2cm}
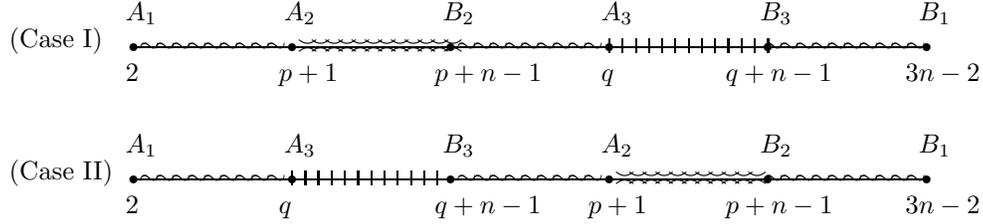
\begin{picture}(0,70)

\put(-17,30){(Case I)}

\put(30,30){\line(1,0){300}}

\put(30,30){\circle*{3}}
\put(90,30){\circle*{3}}
\put(150,30){\circle*{3}}
\put(210,30){\circle*{3}}
\put(270,30){\circle*{3}}
\put(330,30){\circle*{3}}

\put(25,40){$A_1$}
\put(85,40){$A_2$}
\put(145,40){$B_2$}
\put(205,40){$A_3$}
\put(265,40){$B_3$}
\put(325,40){$B_1$}

\put(27,17){2}
\put(83,17){$p+1$}
\put(142,17){$p+n-1$}
\put(205,17){$q$}
\put(252,17){$q+n-1$}
\put(320,17){$3n-2$}

\put(30,28){$\sim$}
\put(36,28){$\sim$}
\put(42,28){$\sim$}
\put(48,28){$\sim$}
\put(54,28){$\sim$}
\put(60,28){$\sim$}
\put(66,28){$\sim$}
\put(72,28){$\sim$}
\put(78,28){$\sim$}
\put(266,28){$\sim$}
\put(272,28){$\sim$}
\put(278,28){$\sim$}
\put(284,28){$\sim$}
\put(290,28){$\sim$}
\put(296,28){$\sim$}
\put(302,28){$\sim$}
\put(308,28){$\sim$}
\put(314,28){$\sim$}
\put(320,28){$\sim$}

\put(150,28){$\sim$}
\put(156,28){$\sim$}
\put(162,28){$\sim$}
\put(168,28){$\sim$}
\put(174,28){$\sim$}
\put(180,28){$\sim$}
\put(186,28){$\sim$}
\put(192,28){$\sim$}
\put(198,28){$\sim$}

\put(90,28){$\asymp$}
\put(95,28){$\asymp$}
\put(100,28){$\asymp$}
\put(105,28){$\asymp$}
\put(110,28){$\asymp$}
\put(115,28){$\asymp$}
\put(120,28){$\asymp$}
\put(125,28){$\asymp$}
\put(130,28){$\asymp$}
\put(135,28){$\asymp$}
\put(140,28){$\asymp$}
\put(145,28){$\asymp$}

\put(210,28){\line(0,1){5}}
\put(215,28){\line(0,1){5}}
\put(220,28){\line(0,1){5}}
\put(225,28){\line(0,1){5}}
\put(230,28){\line(0,1){5}}
\put(235,28){\line(0,1){5}}
\put(240,28){\line(0,1){5}}
\put(245,28){\line(0,1){5}}
\put(250,28){\line(0,1){5}}
\put(255,28){\line(0,1){5}}
\put(260,28){\line(0,1){5}}
\put(265,28){\line(0,1){5}}
\put(270,28){\line(0,1){5}}

\put(-17,-20){(Case II)}

\put(30,-20){\line(1,0){300}}

\put(30,-20){\circle*{3}}
\put(90,-20){\circle*{3}}
\put(150,-20){\circle*{3}}
\put(210,-20){\circle*{3}}
\put(270,-20){\circle*{3}}
\put(330,-20){\circle*{3}}

\put(25,-10){$A_1$}
\put(85,-10){$A_3$}
\put(145,-10){$B_3$}
\put(205,-10){$A_2$}
\put(265,-10){$B_2$}
\put(325,-10){$B_1$}

\put(27,-33){2}
\put(83,-33){$q$}
\put(142,-33){$q+n-1$}
\put(200,-33){$p+1$}
\put(252,-33){$p+n-1$}
\put(320,-33){$3n-2$}

\put(30,-22){$\sim$}
\put(36,-22){$\sim$}
\put(42,-22){$\sim$}
\put(48,-22){$\sim$}
\put(54,-22){$\sim$}
\put(60,-22){$\sim$}
\put(66,-22){$\sim$}
\put(72,-22){$\sim$}
\put(78,-22){$\sim$}
\put(266,-22){$\sim$}
\put(272,-22){$\sim$}
\put(278,-22){$\sim$}
\put(284,-22){$\sim$}
\put(290,-22){$\sim$}
\put(296,-22){$\sim$}
\put(302,-22){$\sim$}
\put(308,-22){$\sim$}
\put(314,-22){$\sim$}
\put(320,-22){$\sim$}

\put(150,-22){$\sim$}
\put(156,-22){$\sim$}
\put(162,-22){$\sim$}
\put(168,-22){$\sim$}
\put(174,-22){$\sim$}
\put(180,-22){$\sim$}
\put(186,-22){$\sim$}
\put(192,-22){$\sim$}
\put(198,-22){$\sim$}

\put(90,-22){\line(0,1){5}}
\put(95,-22){\line(0,1){5}}
\put(100,-22){\line(0,1){5}}
\put(105,-22){\line(0,1){5}}

\put(110,-22){\line(0,1){5}}
\put(115,-22){\line(0,1){5}}
\put(120,-22){\line(0,1){5}}
\put(125,-22){\line(0,1){5}}
\put(130,-22){\line(0,1){5}}
\put(135,-22){\line(0,1){5}}
\put(140,-22){\line(0,1){5}}
\put(145,-22){\line(0,1){5}}

\put(210,-22){$\asymp$}
\put(215,-22){$\asymp$}
\put(220,-22){$\asymp$}
\put(225,-22){$\asymp$}
\put(230,-22){$\asymp$}
\put(235,-22){$\asymp$}
\put(240,-22){$\asymp$}
\put(245,-22){$\asymp$}
\put(250,-22){$\asymp$}
\put(255,-22){$\asymp$}
\put(260,-22){$\asymp$}
\end{picture}
\end{minipage}
$

\vspace{2cm}

Note that 
$|X_1|=n-2, |X_2|=n, |X_3|=n.$

In Case I we have an element 
$$[e](X_1,X_2,X_3)=$$
$$[\omega](t_2,\cdots,t_p,t_{p+n-1},\ldots,t_{q-1},t_{q+n},\ldots,t_{3n-2},[\omega](t_1,t_{p+1},\ldots,t_{p+n-1}),[\omega](t_q,\ldots,t_{q+n-1}))$$
with signature 
$$(-1)^{|X_1|+|(B_2,A_3)||[A_2,B_2]|+|(B_3,B_1]||[A_2,B_2]|+|(B_3,B_1]||[A_3,B_3]|}.$$
Note that
$$|X_1|=n-2\equiv n\mod 2,$$
$$|(B_2,A_3)|=|[B_2,A_3]|-2\equiv |[B_2,A_3]|=q-p-n \mod 2,$$
$$|(B_3,B_1]|= |[B_3,B_1]|-1= 3n-2-q+n-1\equiv q+1\mod\,2,$$
$$|[A_2,B_2]|=\equiv n-1\mod\,2,$$
$$|[A_3,B_3]|\equiv n\mod\,2.$$
Therefore, in Case I, $q\ge p+n,$ and
$$sign\, X_1X_2X_3=(-1)^{n+(q-p-n)(n-1)+(q+1)(n-1)+(q+1)n}=(-1)^{(q-p+1)n+p+1}.$$

In Case II $p\ge q+n-1$ and we have an element 
$$[e](X_1,X_2,X_3)=$$
$$[\omega](t_2,\cdots,t_{q-1},t_{q+n},\ldots,t_{p},t_{p+n},\ldots,t_{3n-2},[\omega](t_1,t_{p+1},\ldots,t_{p+n-1}),[\omega](t_q,\ldots,t_{q+n-1}))$$
with signature 
$$(-1)^{|X_1|+|(B_3,A_2)||[A_3,B_3]|+|(B_2,B_1]||[A_3,B_3]|+|(B_2,B_1]||[A_2,B_2]|+|[A_2,B_2]||[A_3,B_3]|}$$
Since
$$|X_1|\equiv n\mod 2,\quad |[A_3,B_3]|=n,\quad |[A_2,B_2]|\equiv n-1\mod 2,$$
$$|(B_3,A_2)|\equiv |[B_3,A_2]|\equiv p-q-n+1\mod 2,$$
$$|(B_2,B_1]|=|[B_2,B_1]|-1\equiv 3n-2-p-n+1\equiv p-1\mod 2,$$
we have 
$$sign\,X_1X_2X_3=(-1)^{n+(p-q-n+1)n+(p-1)n+(p-1)(n-1)+(n-1)n}=(-1)^{(p-q)n+p-1+n}.$$
So, 
\begin{equation}\label{30}
sign\,X_1X_2X_3=\left\{\begin{array}{ll}
(-1)^{(q+p+1)n+p+1}& \mbox{ Case I, $q\ge p+n$}\\
&\\
(-1)^{(p+q+1)n+p+1}&\mbox{ Case II, $p\ge q+n$}\\
\end{array}\right.
\end{equation}

For $1\le i\le 3n-2$  denote by $\mu_{2,i}^{(p,q)},$
 the coefficient at  $e_i=e(2,\ldots,i,1,i+1,\ldots,3n-2)$
of the element  
$$[\omega](t_2,\cdots,t_p,t_{p+n-1},\ldots,t_{q-1},t_{q+n},\ldots,t_{3n-2},[\omega](t_1,t_{p+1},\ldots,t_{p+n-1}),[\omega](t_q,\ldots,t_{q+n-1})),$$ 
in Case I  or of the element  
$$[\omega](t_2,\cdots,t_{q-1},t_{q+n},\ldots,t_{p},t_{p+n-1},\ldots,t_{3n-2},[\omega](t_1,t_{p+1},\ldots,t_{p+n-1}),[\omega](t_q,\ldots,t_{q+n-1})), $$
in Case II.

To calculate $\mu_{2,i}^{(p,q)}e_i$  we have to do the following permutations.
 
 \bigskip 
 
 Case I. $p+n\le q,$ $p\le i\le p+n-1.$ 

$$[\omega](t_2,\cdots,t_p,t_{p+n},\ldots,t_{q-1},t_{q+n},\ldots,t_{3n-2},[\omega](t_1,t_{p+1},\ldots,t_{p+n-1}),[\omega](t_q,\ldots,t_{q+n-1}))\rightsquigarrow$$

$$
(-1)^{p+n+1}
[\omega](t_2,\cdots,t_p,
[\omega](t_1,t_{p+1},\ldots,t_{p+n-1}),
t_{p+n},\ldots,t_{q-1},t_{q+n},\ldots,t_{3n-2},
\omega(t_q,\ldots,t_{q+n-1}))\rightsquigarrow$$

$$
(-1)^{p+q+n}
[\omega](t_2,\cdots,t_p,
[\omega](t_1,t_{p+1},\ldots,t_{p+n-1}),
t_{p+n},\ldots,t_{q-1},
\omega(t_q,\ldots,t_{q+n-1}),
t_{q+n},\ldots,t_{3n-2})
\rightsquigarrow
$$

$$
(-1)^{i+q+n}
\omega(t_2,\cdots,t_p,
\omega(t_{p+1},\ldots,t_i,t_1,t_{i+1},\ldots, t_{p+n-1}),
t_{p+n},\ldots,t_{q-1},$$
$\qquad\qquad\qquad\qquad\qquad\qquad\qquad\qquad\qquad\qquad\qquad\qquad
\omega(t_q,\ldots,t_{q+n-1}),
t_{q+n},\ldots,t_{3n-2})\rightsquigarrow
$

(by total associativeity)

$$\rightsquigarrow (-1)^{i+q+n}e(2,\ldots,i,1,i+1,\ldots,3n-2)=(-1)^{i+q+n}e_i.$$

\bigskip

Case II, $q+n-1\le p,$ $ p\le i\le p+n-1.$

$$[\omega](t_2,\cdots,t_{q-1},t_{q+n},\ldots,t_{p},t_{p+n},\ldots,t_{3n-2},[\omega](t_1,t_{p+1},\ldots,t_{p+n-1}),[\omega](t_q,\ldots,t_{q+n-1}))\rightsquigarrow$$

$$
(-1)^{p+1}
[\omega](t_2,\cdots,t_{q-1},t_{q+n},\ldots,t_{p},
[\omega](t_1,t_{p+1},\ldots,t_{p+n-1}),
t_{p+n},\ldots,t_{3n-2},
\omega(t_q,\ldots,t_{q+n-1}))\rightsquigarrow$$

$$
(-1)^{p-q-n}
[\omega](t_2,\cdots,t_{q-1},
\omega(t_q,\ldots,t_{q+n-1}),
t_{q+n},\ldots,t_{p},
[\omega](t_1,t_{p+1},\ldots,t_{p+n-1}),
t_{p+n},\ldots,t_{3n-2})\rightsquigarrow$$

$$
(-1)^{i-q-n}
\omega(t_2,\cdots,t_{q-1},
\omega(t_q,\ldots,t_{q+n-1}),t_{q+n},\ldots,t_{p},
\omega(t_{p+1},\ldots,t_i,t_1,t_{i+1},\ldots, t_{p+n-1}),$$
$\hspace{11cm}
t_{p+n},\ldots,t_{3n-2})\rightsquigarrow$

(by total associativeity)

$$\rightsquigarrow (-1)^{i-q-n}e(2,\ldots,i,1,i+1,\ldots,3n-2)=(-1)^{i-q-n}e_i.$$

\medskip

\sudda{
Note that $\mu_{2,i}^{(p,q)}=0$ if $i\notin X_2.$ Therefore,   
$$\mu_{2,i}^{(p,q)}=\left\{ \begin{array}{ll}
0& \mbox{ if } i+1\le p\le q-n \mbox{ or } p\le Min \{i-n+1,q+n-1\}\\
(-1^{i+q+n} & \mbox{ if } p\le Min\{i,q+n-1\} \mbox{ and } i\le p+n-1
\end{array}\right.
$$
if $p\le q-n<q$ 
and 
$$\mu_{2,i}^{(p,q)}=\left\{ \begin{array}{ll}
0& \mbox{ if } i+1\le p\le q+n-1 \mbox{ or } p\le Min \{i-n+1,q-n\}\\
(-1^{i+q+n} & \mbox{ if } p\le Min\{i,q-n\} \mbox{ and } i\le p+n-1
\end{array}\right.
$$
if $q\le p-n+1$
}

Consider the case $i\le n.$ Then the Case II is impossible,  $P_2$ is on the left of $P_3.$ 
We have 
$$\mu_{2,i}=\sum_{p,q} \mu_{2,i}^{(p,q)}= \sum_{p=1}^{i}\sum_{q=p+n}^{2n-1}(-1)^{(q+p+1)n+p+1}(-1)^{i+q+n}=
\sum_{p=1}^{i}\sum_{q=p+n}^{2n-1}(-1)^{(q+p+1)(n+1)+i+n}$$ 
So, if $n$ is even,
$$\mu_{2,i}=
\sum_{p=1}^{i}\sum_{q=p+n}^{2n-1}(-1)^{q+p+1+i}=
(-1)^{i+1}\sum_{p=1}^{i}\sum_{q=p+n}^{2n-1}(-1)^{q-p}=
(-1)^{i + 1}\lfloor \frac{i + 1}{2}\rfloor,$$
and  if $n$ is odd,
$$\mu_{2,i}=
\sum_{p=1}^{i}\sum_{q=p+n}^{2n-1}(-1)^{i+1}=
(-1)^{i + 1}\, \frac{(2 n - i - 1) i}{2}
.$$

Consider the case $n+1\le i\le 2n-2.$ 
Then 
$$\mu_{2,i}=\sum_{p,q} \mu_{2,i}^{(p,q)}= 
\sum_{p=1}^i\sum_{q=p+n}^{2n-1} (-1)^{(q+p+1)n+p+1}(-1)^{i+q+n}+\sum_{p=i-n+1}^{i}\sum_{q=2}^{p-n+1}
(-1)^{(q+p+1)n+p+1}(-1)^{i-q-n}=
$$
$$
(-1)^i(\sum_{p=1}^i\sum_{q=p+n}^{2n-1} (-1)^{(q+p+1)(n+1)+n}+\sum_{p=i-n+1}^{i}\sum_{q=2}^{p-n+1}
(-1)^{(q+p+1)(n+1)+n})
$$
So, if $n$ is even, then 
$$\mu_{2,i} = 
(-1)^i(\sum_{p=i-n+1}^i\sum_{q=p+n}^{2n-1} (-1)^{q+p+1}+\sum_{p=i-n+1}^{i}\sum_{q=2}^{p-n+1}
(-1)^{q+p+1}) =
(-1)^{i+1}(\frac{n}{2}+2\lfloor\frac{n - i-1}{2}\rfloor),$$
and if $n$ is odd 
$$\mu_{2,i}=(-1)^{i+1}(\sum_{p=i-n+1}^i\sum_{q=p+n}^{2n-1} 1+
\sum_{p=i-n+1}^{i}\sum_{q=2}^{p-n+1}
1) =  (-1)^{i + 1} (\frac{n (n - 1)}{2} + (i - n) (-2 n + i + 1)).$$

Consider the case $2n-1\le i\le 3n-2.$ Then Case I is impossible, and $P_2$ is on the right sight  of $P_3.$ We have 
$$\mu_{2,i}=\sum_{p,q} \mu_{2,i}^{(p,q)}=  
\sum_{p=i-n+1}^{2n-1}\sum_{q=2}^{p-n+1}
(-1)^{(q+p+1)n+p+1}(-1)^{i-q-n}=
\sum_{p=i-n+1}^{2n-1}\sum_{q=2}^{p-n+1}
(-1)^{(q+p+1)(n+1)+i+n}
$$
So, if $n$ is even,
$$\mu_{2,i}=(-1)^{i+1}
\sum_{p=i-n+1}^{2n-1}\sum_{q=2}^{p-n+1}
(-1)^{q-p}=
(-1)^i \lfloor \frac{3 n - i}{2}\rfloor$$
and  if $n$ is odd,
$$\mu_{2,i}=(-1)^{i+1}\sum_{p=i-n+1}^{2n-1} \sum_{q=2}^{p-n+1}1=
(-1)^i \frac{ (3n - i - 1) (n - i)}{2}.$$
Lemma \ref{41} is proved completely. 

{\bf Proof of Theorem   \ref{2}}. It follows from Lemmas \ref{40} and \ref{41}.

\section{\bf Proof of Theorem \ref{3}}

Repeats arguments of the proof of Theorem \ref{2}.  
Let 
$$\mu^+_{i}=\left\{ \begin{array}{ll}
\frac{(2n-i-1)i}{2} & \mbox{ if $i\le n$ }\\
&\\
 \frac{n(n-1)}{2}+(i-n)(-2n+i+1) & \mbox{ if $ n+1\le  i\le 2n-2$ }\\
&\\
 \frac{(3n-i-1)(i-n)}{2}& \mbox{ if $2n-1\le i\le 3n-2$ }\\
\end{array}\right.
$$

\begin{lm}\label{50} 
Let $\mu_{1,i}^+$ be the coefficient at $e_i$ of the element 
$F^{\{2\}}_1(t_1,\ldots,t_{3n-2}).$ Then 
$$\mu_{1,i}^+=\mu_i^+,$$
for any $1\le i\le 3n-2.$ 
\end{lm}

\begin{lm}\label{51} 
Let $\mu_{2,i}^+$ be the coefficient at $e_i$ of the element 
$F^{\{2\}}_2(t_1,\ldots,t_{3n-2}).$ Then 
$$\mu_{2,i}^+=\mu_i^+,$$
for any $1\le i\le 3n-2.$ 
\end{lm}

Theorem \ref{3} follows from Lemmas \ref{50} and \ref{51}.

{\bf Remark.} Note that $\mu_i=\mu_i^+(-1)^{i+1},$ $1\le i\le 3n-2,$ if $n$ is odd.   
The  generating function for $\mu_i^+$ is a product of two polynomials,
$$G_n(x)=\sum_{i=1}^{3n-2} \mu_i^+x^{i}=(\sum_{i=1}^n x^i)(\sum_{i=1}^{n-1} (n-i)x^{i-1}+i \,x^{i+n-1}).$$
or,
$$G_n(x)=\sum_{i=1}^{3n-2} \mu_i^+x^{i}=(\sum_{i=1}^n x^i)(n x^{-1}+(x^n-1)\der)(\sum_{i=1}^{n-1}x^i).$$

 If  $n=2k$ is even,  then the generating function for $\mu_i$ is the following polynomial 
$$Q_{2k}(x)=\sum_{i=1}^{6k-2} \mu_{i}x^i=(x-1)^2x(x+1)(\sum_{i=1}^{k} x^{2i-2})^3,$$
or,
$$Q_n(x)=\sum_{i=1}^{3n-2} \mu_{i}x^i=\frac{x(x-1)(x^{n}-1)^3}{(1-x^2)^2}.$$
Therefore, we can formulate the following  more exact versions of Theorems \ref{2}, \ref{3}.
$$F^{[2]}_1=F^{[2]}_2=\sum_{i=1}^{3n-2} \mu_i[e_i],$$
$$F^{\{2\}}_1=F^{\{2\}}_2=\sum_{i=1}^{3n-2} \mu_i^+\{e_i\},$$
 where $\mu_i^+$  are defined as coefficients of the polynomial $G_n(x)$ 
and $\mu_i$ are coefficients of the polynomial $-G_n(-x)$
for odd $n$ and  $Q_n(x)$ for even $n.$

\end{document}